\documentclass[12pt]{amsart}
\usepackage{amscd,amssymb}
\usepackage[arrow,matrix]{xy}
\usepackage[plainpages,backref,urlcolor=blue]{hyperref}

\theoremstyle{plain}
\newtheorem{thm}[subsection]{Theorem}
\newtheorem{lem}[subsection]{Lemma}
\newtheorem{prop}[subsection]{Proposition}
\newtheorem{cor}[subsection]{Corollary}

\theoremstyle{definition}
\newtheorem{rk}[subsection]{Remark}
\newtheorem{definition}[subsection]{Definition}
\newtheorem{ex}[subsection]{Example}
\newtheorem{question}[subsection]{Question}

\numberwithin{equation}{section}
\newcommand{\OO}{{\mathcal O}}
\newcommand{\J}{{\mathcal J}}

\newcommand{\A}{{\mathcal A}}

\newcommand{\wJ}{\widehat{J}}

\newcommand{\CC}{{\mathcal C}}

\newcommand{\V}{{\mathcal V}}

\newcommand{\NN}{{\mathcal N}}

\newcommand{\Z}{\mathbb{Z}}

\newcommand{\C}{\mathbb{C}}

\newcommand{\PP}{\mathbb{P}}

\DeclareMathOperator{\defect}{def}

\begin{document}
\date{20.01.2014 }

\title [Syzygies and logarithmic vector fields along plane curves]
{Syzygies and logarithmic vector fields along plane curves  }

\author[Alexandru Dimca]{Alexandru Dimca$^1$}
\address{Univ. Nice Sophia Antipolis, CNRS,  LJAD, UMR 7351, 06100 Nice, France.}
\email{dimca@unice.fr}

\author[Edoardo Sernesi]{Edoardo Sernesi$^2$}
\address{Dipartimento di Matematica e Fisica, Universit\`a Roma Tre, Largo S. L. Murialdo 1, 00146 Roma, Italy.}
\email{sernesi@mat.uniroma3.it }

\thanks{$^1$ Partially supported by Institut Universitaire de France.} 
\thanks{$^2$ Partially supported by  Project MIUR-PRIN 2010/11 \emph{Geometria delle variet\`a algebriche}}

\subjclass[2000]{Primary 14C34, Secondary 14H50, 32S05 }

\keywords{syzygy, plane curve, logarithmic vector fields, stable bundle, free divisor, Torelli property }

\begin{abstract} We investigate the relations between the syzygies of the Jacobian ideal  of the defining equation for a plane curve $C$ and the stability of the sheaf of logarithmic vector fields along $C$, the freeness of the divisor $C$ and the Torelli properties  of $C$ (in the sense of Dolgachev-Kapranov). We show in particular that curves with a small number of nodes and cusps are Torelli in this sense.

\end{abstract}

\maketitle

%\tableofcontents

\section{Introduction} \label{sec:intro}
Let $C:f=0$ be a complex  projective plane curve, having only weighted homogeneous singularities.  In this paper we continue the investigation of the relations between the syzygies of the Jacobian ideal $J_f$ of $f$ and the stability of the sheaf of logarithmic vector fields $T\langle C\rangle=Der(-logC)$ along $C$, the freeness of the divisor $C$ and Torelli properties  of $C$ (in the sense of Dolgachev-Kapranov \cite{DK})  started by the second author in \cite{Se}.

In the second section we state for reader's convenience as Theorem \ref{vanishing} a result from 
\cite{DS}, giving a sharp lower bound for the degree of (homogeneous) syzygies among the partial derivatives $f_x,f_y,f_z$ of the polynomial $f$ in terms of the Arnold exponents of the singular points $p$ of $C$. Some consequences on the position of singularities of $C$, expressed in terms of defects of linear systems, are also given.

In the third section we recall the definition and basic properties of the sheaf of logarithmic vector fields $T\langle C\rangle$ along $C$, which is in fact a rank two vector bundle on $\PP^2$ in this case. 
For a nodal curve having irreducible components $C_1,\dots,C_r$ whose normalizations are $\widetilde{C}_1,\dots,\widetilde{C}_r$, we prove the formula 
 \[
h^1(T\langle C\rangle(-3)) = h^1(T\langle C\rangle(d-3)) = \sum_i g(\widetilde{C}_i),
\]
where $g(\widetilde{C}_i)$ denotes the genus of $\widetilde{C}_i$, see Proposition \ref{genus}.

In the next section,  we obtain an easy to check sufficient condition for  the stability of the vector bundle $T\langle C\rangle$ expressed as an inequality involving the degree $d$ of $C$ (supposed to have only simple singularities) and the above Arnold exponents, see  Theorem \ref{stabthm} .
Then  we derive a consequence of Theorem \ref{stabthm}, see Corollary  \ref{freecor},
which shows that a curve with a given list of simple singularities is not free, i.e. $T\langle C\rangle$ is not the direct sum of two line bundles, if the degree $d$ of $C$ is large enough. For reader's convenience we include in Remark \ref{rkalgebra} a discussion on the algebraic vs. geometric approaches to the freeness of a plane curve.
To end this section, we discuss two examples, the first one  in common with \cite{ST} and \cite{BC}, of families of curves (with degrees as large as we like)
which are neither free nor stable.

The last section is devoted to Torelli-type questions. After the definition of a Torelli-type curve
(in the sense of Dolgachev-Kapranov), we show that the natural map from the Severi variety of plane reduced curves with a fixed number of nodes $n$ and of cusps $\kappa$ to the corresponding moduli space of stable rank $2$ vector bundles on $\PP^2$ is a morphism, see Proposition \ref{severi}.
This allows us to reprove the known fact that certain reduced curves with many nodes are not Torelli.

On the other side, we {\it conjecture  that any irreducible nodal curve is Torelli.} If the curve is smooth
(and not of Sebastiani-Thom type), this result was established by Ueda and Yoshinaga in \cite{UY}
(where the smooth hypersurface case is treated).
The main result of this paper says that the above conjecture holds for curves with a small number of nodes, i.e. if $n \le (d-1)/2$, see Theorem \ref{Torthm}. A more precise statement is given in 
Theorem \ref{nodalthm}. A version covering curves with few nodes and cusps is given in 
Theorem \ref{nodesandcuspsthm}. We note that  irreducible cuspidal curves are not Torelli in general. The explicit example of a sextic with nine cusps is discussed in detail.

\section{A vanishing result for syzygies among $f_x,f_y$ and $f_z$} \label{sec1}

Let $f$ be a homogeneous polynomial of degree $d$ in the polynomial ring $S=\C[x,y,z]$ and denote by $f_x,f_y,f_z$ the corresponding partial derivatives.
One can consider the graded $S-$submodule $AR(f) \subset S^{3}$ of {\it all relations} involving these derivativess, namely
$$\rho=(a,b,c) \in AR(f)_m$$
if and only if  $af_x+bf_y+cf_z=0$ and $a,b,c$ are in $S_m$.
Let $C$ be the plane curve in $\PP^2$ defined by $f=0$ and assume that $C$ is reduced.
Let $\alpha_C$ be the minimum of the Arnold exponents (alias singularity indices or log canonical thresholds, see Theorem 9.5 in \cite{Ko}) $\alpha_p$ of the singular points $p$ of $C$. If the germ $(C,p)$ is weighted homogeneous of type $(w_1,w_2;1)$ with $0<w_j \leq 1/2$, then one has 
\begin{equation}
\label{Ae}
\alpha_p={w_1}+{w_2},
\end{equation}
see for instance \cite{DS0}. Moreover, since for any isolated plane curve singularity $(C,0)$ the spectrum of $(C,p)$ is contained in the interval $(0,2)$ and it is symmetric with respect to $1$, it follows that $\alpha_p \leq 1$ with equality exactly when 
$(C,p)$ is a node, i.e. an $A_1$-singularity.
With this notation, Corollary 5.5 in \cite{DS} can be restated as follows.

\begin{thm}
\label{vanishing}
Let $C:f=0$ be a degree $d$ reduced curve in $\PP^2$ having only weighted homogeneous singularities. 
Then $AR(f)_m=0$ for all $m <\alpha_C d-2$.
\end{thm}

\proof

It is enough to check that one has the obvious identification $N_{d+k}=AR(f)_{k-2}$, for any $k<d+1$, where the graded $S$-module $N$ is defined in \cite{DS} using a shifted version of the Koszul complex for $f_x,f_y,f_z$.

\endproof 

\begin{ex}
\label{exvanishing}
(i) If $C:f=0$ is a degree $d$ nodal curve in $\PP^2$, then
$ \alpha_C =1$.
It follows that $AR(f)_m=0$ for all $m \leq d-3$ which is exactly the bound obtained in Thm. 4.1 in \cite{DSt}. This bound is known to be optimal, since $\dim AR(f)_{d-2}=r-1$, where $r$ is the number of irreducible components of $C$,  see Thm. 4.1 in \cite{DSt}.

(ii) If $C:f=0$ is a degree $d$  curve in $\PP^2$ having only nodes $A_1$ and cusps $A_2$ as singularities, then
$ \alpha_C =5/6$. It follows that $AR(f)_m=0$ for all $m < 5d/6-2$. For the Zariski sextic curve with 6 cusps on a conic, e.g. $f=(x^2+y^2)^3+(y^3+z^3)^2$, this bound is sharp since $AR(f)_3\ne 0$. As in (i) above, such non vanishing results have a geometrical meaning (at least in many cases). For instance $AR(f)_3\ne 0$ in the case of the Zariski sextic is related to the fact that the action of the monodromy on $H^1(F,\C)$ is not the identity, where $F:f(x,y,z)=1$ denotes the Milnor fiber of the defining equation $f$, see \cite{DS} for the general theory. Similar remarks apply to the non-vanishing claimed in the following point (iii).

(iii) If $C:f=0$ is a degree $d$  curve in $\PP^2$ having only nodes $A_1$, cusps $A_2$  and ordinary triple points $D_4$ as singularities, then
$ \alpha_C =2/3$. It follows that $AR(f)_m=0$ for all $m < 2d/3-2$. For the line arrangement defined by $f=(x^2-y^2)(y^2-z^2)(z^2-x^2)$ and for the curve $f=(x^3+y^3+z^3)^3+(x^3+2y^3+3z^3)^3$ with 3 irreducible components (each smooth of genus 1) and nine $D_4$ singularities,
this bound is  sharp.

\end{ex}

The following result is also useful in the sequel.

\begin{prop}
\label{simple}
Let $C:f=0$ be a reduced curve in $\PP^2$ having only weighted homogeneous singularities. 
Then  $\alpha_C >1/2$ if and only if $C$ has only simple singularities, i.e. singularities of type $A_k$ for $k \geq 1$, $D_k$ for $k\geq 4$, $E_6$, $E_7$ and $E_8$.
\end{prop}

\proof
Using formula \eqref{Ae}, this is a classical result in singularity theory, see \cite{KS}. One can also look at Corollary 7.45 and its proof in \cite{D0}.

\endproof 

\begin{rk}
\label{rkalpha}

 For a non weighted homogeneous plane curve singularity $(C,p)$, the computation of the corresponding exponent $\alpha_p$ is much more complicated. For instance recall that we have  $\alpha_p= 1/2$ for any singularity $(C,p)$ in the series of unimodal singularities $T_{2,q,r}: x^q+x^2y^2+y^r=0$, where $q \geq 2$, $r\geq 2$, see \cite{AGV}, Table 2, page 189.

\end{rk}

We discuss now some consequences of the above results on 
the position of the singularities of $C$, which is described to some extent by the sequence of defects
\begin{equation}
\label{defect}
\defect _k\Sigma_f=\tau(C)-\dim \frac{S_k}{\wJ_k}.
\end{equation}
Here $\Sigma_f$ is the subscheme of $\PP^2$ defined by the Jacobian ideal $J=J_f$  and $\wJ=\wJ_f$ denotes the saturation of $J$. Set $R(f)=S/J$, the corresponding graded Jacobian (or Milnor) algebra of $f$, see \cite{D2} and note that 
\[
\frac{S_k}{\wJ_k} = H^1_\mathbf{m}(R(f))_k.
\]
This  follows from the exact sequence (11), p. 7, of \cite{Se}, because 
$$\frac{S}{\wJ}=R(f)/H^0_\mathbf{m}(R(f)).$$
Other related invariants of the curve $C$ have been introduced in \cite{DSt}, and we recall them below.

\begin{definition}
\label{def}
Let $C:f=0$ a degree $d$ curve  with isolated singularities in $\PP^2$.

\noindent (i) the {\it coincidence threshold} $ct(C)$  is defined as
$$ct(C)=\max \{q~~:~~\dim M(f)_k=\dim M(f_s)_k \text{ for all } k \leq q\},$$
with $f_s$  a homogeneous polynomial in $S$ of degree $d$ such that $C_s:f_s=0$ is a smooth curve in $\PP^2$.

\noindent (ii) the {\it minimal degree of a nontrivial relation} $mdr(D)$ is defined as
$$mdr(C)=\min \{q~~:~~ ER(f)_q\ne 0\},$$
where $ER(f)$ is the quotient of the graded $S$-module $AR(f)$ by the submodule spanned by the Koszul (trivial) relations among $f_x,f_y,f_z$.

\end{definition}
It is known that one has
\begin{equation} 
\label{REL}
ct(C)=mdr(C)+d-2,
\end{equation} 
see \cite{DSt}, formula (1.3).

\begin{ex}
\label{exmdr}
 If $C:f=0$ is a degree $d$ nodal curve in $\PP^2$, then 
$$d-2 \leq mdr(C) \leq 2(d-2) \text { and }  2(d-2) \leq ct(C) \leq 3(d-2) .$$
Moreover the equalities $d-2 = mdr(C)$, $2(d-2)=ct(C)$ hold exactly when $C$ is not irreducible, and the equalities
$mdr(C) = 2d-4$, $ct(C)=3(d-2)$ hold exactly when $C$ has just one node, see \cite{DSt}.

\end{ex}

More generally, the following result holds.
\begin{prop}
\label{position}
Let $C:f=0$ be a reduced curve in $\PP^2$ having only weighted homogeneous singularities. 
Then the following holds.

\medskip

\noindent (i) $ct(C) \geq (\alpha_C+1)d-4$; in particular, if $C$ has only simple singularities, then $ct(C) \geq 3d/2-3=T/2$, with $T=3(d-2)=\max \{q~:~ M(f_s)_q \ne 0\}$.

\medskip

\noindent (ii) $\defect_k\Sigma_f=0$ for $k \geq (2-\alpha_C)d-2$; in particular, if $C$ has only simple singularities, then  $\defect_k\Sigma_f=0$ for   $k \geq 3d/2-2=T/2+1$.
\end{prop}

\proof
The first claim is a direct consequence of formula \eqref{REL} and Theorem 
\ref{vanishing}. The second claim follows from Theorem 1 in  \cite{D2}, which implies that
$\defect_k\Sigma_f=0$ for $T-k \leq ct(C)$. 
\endproof 

\begin{ex}
\label{exdefect}
(i) If $C:f=0$ is a degree $d$ nodal curve in $\PP^2$, then $ct(C) \geq 2d-4$ and $\defect_k\Sigma_f=0$ for $k \geq d-2$ were already obtained in
  \cite{DSt} and are sharp. 
Actually  the nodes of an irreducible plane curve $C$ impose independent conditions to the curves of degree $k\ge d-3$ and therefore 
	$\defect_k\Sigma_f=0$ for $k \ge d-3$ in this case. This a consequence of the theorem of Gorenstein (\cite{lS79}, p. 38) and of the fact that in the nodal case the Tjurina, Milnor and adjoint ideals coincide. See also \cite{ACGH}, Ex. 11 p. 54.
	 An example for which $\defect_{d-4}\Sigma_f\ne 0$ is a sextic curve $C$ of genus $4$  which is the projection  of a 
	canonical sextic in $\PP^3$. This curve $C$ has $6$ nodes situated on a conic (see \cite{ACGH}, Ex. 24 p. 57).

(ii) If $C:f=0$ is a degree $d$  curve in $\PP^2$ having only nodes $A_1$ and cusps $A_2$ as singularities,  then $ct(C) \geq 11d/6-4$ and $\defect_k\Sigma_f=0$ for $k \geq 7d/6-2$.

(iii) If $C:f=0$ is a degree $d$  curve in $\PP^2$ having only nodes $A_1$, cusps $A_2$  and ordinary triple points $D_4$ as singularities, then
$ct(C) \geq 7d/4-4$ and $\defect_k\Sigma_f=0$ for $k \geq 5d/4-2$.

\end{ex}

%%%%%%%%%%%%%%%%%%%%%%%%%%%%%%%%%%%%%%%%%%%%%%%%%%%%%%%%%%

\section{Syzygies and logarithmic vector fields}

For a reduced projective plane curve $C:f=0$ of degree $d$, let $T\langle C\rangle=Der(-logC)$ denote the sheaf of logarithmic vector fields along $C$. This sheaf, which can be defined more generally for any hypersurface $D$ in $\PP^n$, is always reflexive, see \cite{KS2}.
Moreover, any reflexive sheaf on a smooth surface is free, see \cite{OSS}, Lemma 1.1.10, page 149. Hence in our setting $T\langle C\rangle$ is a rank $2$ vector bundle.

One has the exact sequence:
\[
\xymatrix{0\ar[r] &T\langle C\rangle \ar[r] & \OO_{\PP^2}(1)^3 \ar[r] & \J_f(d) \ar[r] & 0}\]
where $\J_f \subset \OO_{\PP^2}$ is the gradient ideal sheaf of $f$. This gives an identification:
\begin{equation}\label{logseq1}
AR(f)_m = H^0(\PP^2,T\langle C\rangle(m-1))
\end{equation}
for all $m$. 

%Let $\J_{f/C}=\J_f/\I_C \subset \OO_C$ be the jacobian ideal sheaf of $C$ and $T^1_C = \OO_C/\J_{f/C}$ the first cotangent sheaf of $C$. Then, letting $\Omega^1(log C)=T\langle C\rangle^\vee$ be the sheaf of logaritmic differentials along $C$, we have   exact sequences (\cite{iD07}, (2.2) and Prop. 2.1):
%\[\xymatrix{0 \ar[r] & \Omega^1_{\PP^2} \ar[r] & \Omega^1(log C) \ar[r] & Ext^1_{\PP^2}(\J_{f/C}(d),\OO_{\PP^2}) \ar[r] & 0}\]
%and
%\[\xymatrix{0 \ar[r] & \OO_C \ar[r] & Ext^1_{\PP^2}(\J_{f/C}(d),\OO_{\PP^2}) \ar[r] & T^1_C \ar[r] &0}\]
The Chern classes of $T\langle C\rangle(k)$ are:
\begin{equation}\label{Chern}
c_1(T\langle C\rangle(k)) = 3-d+2k, \quad c_2(T\langle
 C\rangle(k)) = d^2-(3+k)d+3+3k+k^2- \tau(C)
\end{equation}
where       $\tau(C) = h^0(T^1_C)$ is the global Tjurina number of $C$.  Moreover one easily computes that:
\[
\chi(T\langle C\rangle(k)) = 3 {k+3 \choose 2} - {d+k+2 \choose 2} + \tau(C).
\]
In the case $k=d-3$ we obtain:
\[
\chi(T\langle C\rangle(d-3)) = - \left[{d-1 \choose 2}- \tau(C)\right].
\]
Moreover, using Serre duality and the identity 
\begin{equation}\label{logforms}
\Omega^1(log C) = T\langle C\rangle (d-3), 
\end{equation}
which follows from Lemma 4.1 in \cite{Se}, we obtain:
\[
h^2(T\langle C\rangle(d-3)) = h^0(\Omega^1(log C))(-d)= h^0(T\langle C\rangle(-3)) =0.
\]
In conclusion we have:
\begin{align}
\dim AR(f)_{d-2} &= h^0(T\langle C\rangle(d-3))  \label{logseq2}\\
&=h^1(T\langle C\rangle(d-3))+\chi(T\langle C\rangle(d-3)) \notag \\
&=h^1(T\langle C\rangle(d-3))- \left[{d-1 \choose 2}- \tau(C)\right].\notag
\end{align}
A similar computation gives:
\begin{align}
\dim AR(f)_{d-3} &= h^0(T\langle C\rangle(d-4))  \label{logseq3}\\
&= h^1(T\langle C\rangle(d-4))- \left[{d \choose 2} - \tau(C) \right]. \notag
\end{align}
The following proposition generalizes to all nodal curves the dimension estimate of Corollary 5.2 of \cite{Se}.

\begin{prop}\label{genus}
If $C$ has only nodes then
\[
h^1(T\langle C\rangle(-3)) = h^1(T\langle C\rangle(d-3)) = \sum_i g(\widetilde{C}_i),
\]
where $C_1, \dots, C_r$ are the irreducible components of $C$ and $\widetilde{C}_i$ is the normalization of $C_i$, $i=1, \dots, r$.
\end{prop} 

\proof  The first equality is a consequence of the self-duality of $H^1_*(T\langle C \rangle)$. Therefore it suffices to prove the second equality.
We have that  $\tau(C)=\delta$,   the number of nodes of $C$. The geometric genus of $C$ is 
\[
g(C)= {d-1 \choose 2}- \delta = \sum g(\widetilde{C}_i)-r+1.
\]
Therefore we need to prove that 
\[
h^1(T\langle C\rangle(d-3)) = {d-1 \choose 2}- \delta+r-1.
\]
Recalling Example \ref{exvanishing}(i), we see that this follows from (\ref{logseq2}).  \endproof

\begin{rk}
\label{rkotherproof}
An alternative proof of Proposition \ref{genus} can be obtained by using the formula 
\ref{logforms} to pass to logarithmic 1-forms and Proposition 4.1 in \cite{DSt0} alongside with basic facts on mixed Hodge theory.

\end{rk}

%%%%%%%%%%%%%%%%%%%%%%%%%%%%%%%%%%%%%%%%%%%%%%%%%%%%%%%%%%%%%%%

\section{Stability of the bundle $T\langle C\rangle$ and freeness of the divisor $C$} 

Recall that for a rank 2 torsion free coherent sheaf  $E$ on the projective space $\PP^n$ the notions of Mumford-Takemoto stability, Gieseker-Maruyama stability and simplicity (i.e. $End(E)=\C$) coincide, see 
\cite{OSS}, and play a key role in the understanding such sheaves. Since $T\langle C\rangle$ is in a natural way a sub bundle of the tangent bundle $T\PP^2$ see \cite{Se}, which is stable, see \cite{OSS}, Thm. 1.3.2, p. 182, it is natural to ask about its stability properties.
One computes easily that   the discriminant of $T\langle C\rangle$ is:
	\[
	\Delta(T\langle C\rangle) = 3(-d^2+2d-1)+4\tau(C)=4\tau(C)-3(d-1)^2
	\]
and $\Delta(T\langle C\rangle)<0$ is a necessary condition for the stability of $T\langle C\rangle$ (\cite{OSS}, p. 168). This condition already puts some restrictions on $\tau(C)$.
 The second author has shown in \cite{Se}, Proposition 2.4 that for a reduced plane curve $C:f=0$ of degree $d $, the torsion free coherent sheaf  $T\langle C\rangle$ is stable  if and only if
\begin{equation}
\label{stab}
AR(f)_m=0 \text{ for } m \leq (d-1)/2.
\end{equation}
This result combined with Theorem \ref{vanishing} and Proposition \ref{simple} yields the following result.

\begin{thm}
\label{stabthm}
Let $C:f=0$ be a  reduced curve in $\PP^2$ of degree $d$ having only simple singularities. 
Then $T\langle C\rangle$ is stable if either $d$ is odd and 
$$d > \frac{3}{2\alpha_C-1},$$
 or
$d$ is even and 
$$d > \frac{2}{2\alpha_C-1}.$$
\end{thm}

\begin{ex}
\label{exstab}
(i) If $C:f=0$ is a degree $d$ nodal curve in $\PP^2$, then
$ \alpha_C =1$ and the above Theorem tells us that $T\langle C\rangle$ is stable if either $d$ is odd and 
$d \geq 5$ 
 or
$d$ is even and $d \geq 4$, i.e. $T\langle C\rangle$ is stable for all $d\geq 4$.

Note that for  $C:f=X_0X_1X_2=0$, the sheaf $T\langle C \rangle$ splits and therefore it is not stable. Hence our result is sharp in this case. On the other hand, Example \ref{exvanishing} (i) and the formula \eqref{stab} show that $T\langle C \rangle$  is stable for a nodal cubic $C$.

(ii) If $C:f=0$ is a degree $d$  curve in $\PP^2$ having only nodes $(A_1)$ and cusps $(A_2)$ as singularities, then
$ \alpha_C =5/6$ and the above Theorem tells us that $T\langle C\rangle$ is stable if either $d$ is odd and 
$d \geq 5$ 
 or
$d$ is even and $d \geq 4$. The only case not covered is $d=3$, $C$ a cuspidal cubic.
 The cuspidal cubic $X_0X_1^2-X_2^3$ has $h^0(T\langle C\rangle)=AR_1 \ne 0$ so $T\langle C \rangle$ is not stable by Lemma 1.2.5 p. 165 in \cite{OSS}. Hence our result is sharp in this case as well.

\end{ex}

%%%%%%%%%%%%%%%%%%%%%%%%%%%%%%%%%%%%%%%%%%%%%%%%%%%%%%%%%%%%%%%

The reduced plane curve $C$ is free (as a divisor) if the vector bundle $T\langle C\rangle$ splits as a direct sum of two line bundles. The formula \eqref{Chern} for the Chern classes of $T\langle C\rangle(-1)$ implies 
the following
result.
\begin{lem}
\label{lem1}
Suppose the curve $C$ is free, and
$$T\langle C\rangle(-1)=\OO(-a)\oplus \OO(-b).$$
Then the integers $a$ and $b$  above are positive and satisfy the system of equations 
$$a+b=d-1, \   \   ab= (d-1)^2 -\tau(C).$$
 In particular, the discriminant $\Delta(T\langle C\rangle)= 4\tau(C)-3(d-1)^2$ of the bundle $T\langle C\rangle$ is a perfect square.
Moreover, one has $a>0$ and $b>0$ except when $C$ is a union of lines passing through one point.
\end{lem}

\proof The only claim that needs some explanation is about the (strict) positivity of $a$ and $b$.
Note that $\tau(C) \leq \mu(C)$, where $\mu(C)$ denotes the sum of the Milnor numbers of the singularities of $C$. It is well known that $\mu(C) \leq (d-1)^2$ with equality exactly when $C$ is a pencil of curves. Hence $ab= (d-1)^2 -\tau(C) \geq 0$, with strict inequality when $C$ is not a line pencil.

\endproof

By definition, it is clear that if a reduced curve $C$ is free, then $T\langle C\rangle$ is not stable.  Therefore Theorem \ref{stabthm} implies the following:

%it was shown in \cite{Se}, formula (33), that the following hold
%\begin{equation}\label{free2} a+b=d-1 \text{ and } ab=(d-1)^2-\tau(C),\end{equation}
%where $d$ is the degree of $C$ and $\tau(C)$ is the total Tjurina number of $C$. We have the following result.

\begin{cor}
\label{freecor}
Let $C:f=0$ be a  reduced curve in $\PP^2$ of degree $d$ having only simple singularities. 
Then $C$ is not free if either $d$ is odd and 
$$d > \frac{3}{2\alpha_C-1},$$
 or
$d$ is even and 
$$d > \frac{2}{2\alpha_C-1}.$$

\end{cor}

%\proof
%If we assume that $T\langle C\rangle(-1)=\OO(-a)\oplus \OO(-b)$, then it follows that 
%$T\langle C\rangle(a-1)=\OO \oplus \OO(a-b)$ has nontrivial global sections. Using the definition of the sheaf $T\langle C\rangle$, it follows that 
%$$H^0(\PP^2, T\langle C\rangle(s-1)=AR(f)_s$$
%for any integer $s$, see for instance the proof of Prop. 2.4. in \cite{Se}. Hence, by Theorem 
%\ref{vanishing}, we have $a \geq \alpha_C d-2$. By formula \eqref{free2}, we get
%$$d-1 \geq 2a \geq 2\alpha_C d-4.$$
%This completes the proof.
%\endproof 

\begin{ex}
\label{exfree}
(i) If $C:f=0$ is a degree $d$ nodal curve in $\PP^2$, then
$ \alpha_C =1$ and the above Theorem tells us that $C$ is not free  if $d>3$. The cases not covered are the following.
If $d=2$, then either $C$ is union of two lines ($a=0$, $b=1$) which is free as a subarrangement of type $\A_X$ of the arrangement $\A:XYZ=0$ (a triangle), see Theorem 4.37 in \cite{OT}, or $C$ is smooth, and then  $T\langle C\rangle$  is stable \cite{Se}. If $d=3$, then either $C$ is  a nodal cubic (not free, since $a+b=2$, $ab=3$ has no integer solution), or $C$ is a triangle, which is free, see Example \ref{exstab}, (i), or $C$ is smooth and we conclude as in the case of smooth conics. Hence the only free nodal curves are two lines $XY=0$ and the triangle $XYZ=0$.

(ii) If $C:f=0$ is a degree $d$  curve in $\PP^2$ having only nodes $(A_1)$ and cusps $(A_2)$ as singularities, then
$ \alpha_C =5/6$ and the above Theorem tells us that $C$ is not free if  $d>4$.
The case $d=3$ leads to a non free curve since the system $a+b=2$, $ab=2$ has no integer solution. The case $d=4$ leads again to non free curves. Indeed the system becomes
$a+b=3$ and $ab=9-\tau(C)$. The only possible integer solution may be $a=1$ and $b=2$, hence $\tau(C)=7$. When $C$ is irreducible, the genus of the normalization $\tilde C$ is given by
$3-n-k$, where $n$ is the number of nodes and $\kappa$ the number of cusps of $C$. Since 
$\tau(C)=n+2\kappa$, we see that $\tau(C)=7$ cannot be realized. The case of a reducible curve $C$
(cuspidal cubic plus a secant) is even simpler to handle.
Hence our result is not sharp in this case.

\end{ex}

As a special case of Corollary \ref{freecor}  we have the following.
\begin{cor}
\label{freeline arr}
Assume that $C$ is the union of lines of a line arrangement $\A$ in $\PP^2$ having only double and triple points. If either $d>9$ is odd, or $d>6$  is even, then the curve $C$ (or, equivalently, the line arrangement $\A$) is not free.

\end{cor}

This result is sharp, since it is known that the following two arrangements
$$(x^2-y^2)(y^2-z^2)(x^2-z^2)=0$$
and
$$(x^3-y^3)(y^3-z^3)(x^3-z^3)=0$$
are free.

\begin{rk}
\label{rkalgebra} From a purely algebraic view-point, the reduced plane curve $C:f=0$ is free if and only if the corresponding Jacobian (or gradient) ideal $J_f$ spanned by the partial derivatives
$f_x,f_y,f_z$ of $f$ in $S$ is a perfect ideal, i.e. the Jacobian ring $R(f)=S/J_f$ is Cohen-Macauley.
Equivalently, $R(f)$ has a Hilbert-Burch minimal free resolution of the form
$$0\to S^2 \to S^3 \to S \to R(f) \to 0.$$
This is the case exactly when $\wJ_f= J_f$, where $\wJ_f$ denotes the saturation of the ideal $J_f$, see \cite{ST}, the line after Prop. 1.9. In other words, $C$ is free if and only if $\wJ_f/J_f=H^0 _\mathbf{m}(R(f))=0$, see \cite{Se} or \cite{DS}. Geometrically, this follows from Horrocks' Theorem, see \cite{OSS}, p.39, saying that the bundle $T\langle C\rangle$ splits if and only if $H^1(\PP^2,T\langle C\rangle(k))=0$ for any integer $k$. Then one uses the isomorphism
$H^1(\PP^2,T\langle C\rangle(k))=H^0 _\mathbf{m}(R(f))_{d+k}$, see Prop. 2.1 in  \cite{Se}.

There is also a notion of free divisor in local analytic geometry. The two notions are related as follows: the projective curve $C:f=0$ is free if and only if the divisor germ $(D,0)$ in $\C^3$ is free, where $D$ denotes the cone over the curve $C$, i.e. the surface singularity defined by $f=0$ in $\C^3$. For more on this equivalence we refer to \cite{Y} and the references there.

Moreover, note that a reduced plane curve $C$ has only weighted homogeneous singularities if and only if the ideal $J_f$ is of linear type, see Prop. 1.6 in \cite{ST}. Jacobian ideals of linear type are also considered in the local analytic theory, see for instance 
\cite{Mac}.

In the local analytic version of the theory there is also a notion of linear free divisor, apparently not related to the ideals of linear type. For more on linear free divisors see \cite{GMNS} and the references there.

\end{rk}

\begin{rk}
\label{rkalg2}  In this remark we review briefly the examples of free divisors constructed by Simis and Toh\u aneanu in \cite{ST}. In Prop. 2.2. they construct a sequence of irreducible free divisors
$$C_d: f_d=y^{d-1}z+x^d+x^2y^{d-2}+a_3xy^{d-1}+a_4y^d=0,$$
where $a_3,a_4 \in \C$ and $d\geq 5$.
The curve $C_d$ has a unique singularity  located at the point $p=(0:0:1)$ and given in the local coordinates $(x,y)$ by the equation
$$y^{d-1}+x^d+x^2y^{d-2}+a_3xy^{d-1}+a_4y^d=0.$$
This singularity is weighted homogeneous  only for $d=5$ and even then it is not simple.
For $d>5$ this singularity is semi-weighted homogeneous and belongs to the same $\mu$-constant family as the associated weighted homogeneous singularity
$y^{d-1}z+x^d=0.$
Since the exponent $\alpha_p$ is constant in $\mu$-constant families, we infer that
$\alpha_p=1/(d-1)+1/d$. It follows that $d\alpha_p \to 2$ when $d \to \infty$.

Other examples of free divisors in \cite{ST} are described in Cor. 2.7 and are obtained by the homogeneization with respect to $z$ of a weighted homogeneous polynomial $g$ in $x,y$.
These divisors are not irreducible and also have non simple singularities, coming either from the singularity of $g$ at the origin or from other singularities. For instance, if we start with $g=x^2y+y^d$, a simple singularity of type $D$, then $f=x^2yz^{d-3}+y^d$ and the singularity at $p=(1:0:0)$ is not simple as soon as $d>5$. Corollary 2.10 in \cite{ST} describe a slightly different construction, called the coneing, leading to free divisor having essentially the same properties as those in Cor. 2.10 from our point of view. 

Similar constructions based on weighted homogeneity with respect to two distinct sets of weights are given by Buchweitz and Conca in \cite{BC}, see especially Thm. 3.5 and Thm. 6.1.

In all these examples it seems that 
$\alpha_C \to 0$ when $\deg (C) \to \infty.$

\end{rk}

%%%%%%%%%%%%%%%%%%%%%%%%%%%%%%%%%%%%%%%%%%%%%%%%%%%%%%%%%%%%%%%

%%%%%%%%%%%%%%%%%%%%%%%%%%%%%%%%%%%%%%%%%%%%%%%%%%%%%%%%%%%%%%%%%%%%%%%%%%%%%%%%%%%%%%%%%%%

To end this section,
we discuss two examples of families of such curves $C$ which are neither free nor stable, the first one is in common with \cite{ST} and \cite{BC}.

\begin{ex} (A Thom-Sebastiani type example, with $mdr(C)=1$)
\label{exThomSeb}

Consider the family of curves $C=C_{a,b}:f=x^ay^b +z^d$ for $a+b=d$, $a>0$, $b>0$.
It follows from Proposition 2.11 (i) in \cite{ST} or Thm. 6.1 in \cite{BC} that $C$ is not a free divisor. From the obvious relation $bxf_x-ayf_y=0$ it follows that $mdr(f)=1$ and hence by \eqref{stab} that 
$T\langle C\rangle$ is not stable for $d\geq 3$. When $a\geq 2$, $b\geq 2$, then the curve $C$ has two singular points located at $(1:0:0)$ and $(0:1:0)$ and it follows that 
$$\alpha_C=1/d+\min (1/a,1/b) \leq 3/d.$$
We get $\limsup (\alpha_C \cdot d) \leq 3$ when $d\to \infty$, with equality when $a$ and $b$ are chosen both close to $d/2$.

\end{ex}

We consider next a  non Thom-Sebastiani type family of irreducible curves $C$, with $mdr(C) \to \infty$ when $d \to \infty$. These curves are not obtained by homogenization of a weighted polynomial in two variables as some examples in \cite{ST}.

\begin{prop}
\label{nonTS}
Consider the family of irreducible curves 
$$C=C_{a,b,c}:f=x^ay^bz^c+y^d+z^d=0,$$ 
for $a+b+c=d$, $a>1$, $b>1$, $c>1$. Then the following hold.
\medskip

\noindent (i) $mdr(C)=\min (d-b,d-c).$

\medskip

\noindent (ii) If either $d=2b-1$, $c=2$ and $a=b-3$, for $b\geq 5$, or
 $d=2b-2$, $c=2$ and $a=b-4$, for $b\geq 6$, then the bundle $T\langle C\rangle$ is not stable.

\medskip

\noindent (iii) If $b=c=\left[d/2\right]$, then the bundle $T\langle C\rangle$ is stable and $C$ has a unique singularity $p=(1:0:0)$  with very large Milnor number, namely $$\mu(C,p)=2d(b-1)+1 \geq (d-1)^2 -d.$$
(iv) All the curves $C=C_{a,b,c}$ are not free.
\end{prop}

\proof

The fact that $C$ is irreducible is equivalent to the irreducibility of the affine curve 
$F: x^ay^b +y^d+1=0$. Since the polynomial $g(x,y)=x^ay^b +y^d$ is weighted homogeneous, it follows that all fibers $g^{-1}(s)$ are isomorphic for $s \ne 0$. In particular, $F$ can be regarded as the generic fiber of $g$, and hence it is irreducible as $g$ is clearly a primitive polynomial, i.e. not of the form $h(g_1(x,y))$, with $h \in \C[t]$ polynomial of degree $>1$.

One has the following obvious syzygy of degree $(d-b)$
$$dy^{d-b}f_x+x^{a-1}z^c(bxf_x-ayf_y)=0,$$
and a similar one of degree $d-c$ replacing $y$ by $z$. In order to prove the first claim it is enough to show that there are no relations of strictly lower degree. We assume that $b \geq c$, and show there are no relations of  degree $<d-b$. The derivative $f_y$ contains the monomial
$dy^{d-1}$. In a relation $uf_x+vf_y+wf_z=0$, this monomial can cancel with terms coming from
$uf_x$ or with terms coming from $wf_z$.
In the first case, the factor $v$ must contain a monomial divisible by $v_1=x^{a-1}z^c$, in the second
case a monomial divisible by $v_2=x^az^{c-1}$ or $z^d$. The last case is excluded because we consider only relation with $\deg u =\deg v =\deg w < d-b$.
Since $(a-1)+c=a+(c-1)=d-b-1$, it is enough to show that $\deg v =\deg v_1=\deg v_2$ yields a contradiction. Indeed, these equalities implies that $vf_y$ contain terms whose degree with respect to $y$ is $b-1$, and such terms cannot cancel with terms coming from $uf_x$ and $wf_z$. This proves the first claim (i).

To have $T\langle C\rangle$  not stable, it is enough by the relation \eqref{stab}  and (i) to have $b \geq c$ and $d-b \leq (d-1)/2$ which is equivalent to
$d\leq 2b-1$. This proves the second claim (ii).

To prove (iii), one has to use again  \eqref{stab}  and (i) to show that $T\langle C\rangle$ is stable. The formula for the Milnor number follows from the fact that $(C,p)$ is a Newton non-degenerated singularity, and hence 
$$\mu(C,p)=2A-2d+1$$
where $A$ denotes the area between the coordinate axes and the Newton boundary of the defining equation $y^bz^b +y^d +z^d=0$ for the singularity $(C,p)$.

Finally, to prove (iv), it is enough  to show that  $\wJ_f/J_f=H^0 _\mathbf{m}(R(f)) \ne 0$, as explained in  \ref{rkalgebra}. Note that $x^{a-1}y^bz^c \in J_f$, which implies $y^d \in J_f$ and $z^d \in J_f$. It follows that $y^bz^c \in \wJ_f$. On the other hand it is clear that $y^bz^c \notin J_f$, since $\deg y^bz^c <d-1$.

\endproof

\begin{rk}
\label{rklimit}

In both cases (ii) and (iii) above, one can compute the exponent of the Newton non-degenerated singularity  $(C,p)$ using the distance between the point $(1,1)$ and the Newton boundary of the singularity $(C,p)$, see Thm. 6.4 on page 150 in  \cite{AGV}. This implies $\alpha_C \cdot d \to 3$ when $d\to \infty$.

In view of this and the final comment in Remark \ref{rkalg2}, it would be interesting to find examples of families of curves $C_d$, with $\deg C_d=d$, such that $C_d$ is free (resp. $T\langle C\rangle$  is not stable) and $\alpha_{C_d} >\epsilon$ for all $d$ and some fixed $\epsilon >0$.

\end{rk}

\begin{rk}
\label{ct+st+free}

For a recent interesting result involving the invariant $ct(C)$ of $C$ introduced in section 2 and the freeness of the divisor $C$, see  \cite{St}.

\end{rk}
%%%%%%%%%%%%%%%%%%%%%%%%%%%%%%%%%%%%%%%%%%%%%%%%%%%%%%%%%%%%%%

\section{Torelli-type questions}

We will adopt the following

\begin{definition}\label{deftor}
A reduced hypersurface $X\subset \PP^r$ is called \emph{LC-Torelli} (where LC stands for local cohomology) if it can be reconstructed from the $\mathbb{C}[X_0,\dots,X_r]$-module $H^1_*(T\langle X \rangle)=\oplus_k H^1(\PP^r,T\langle X \rangle(k)) $.  
 We say that $X$ is  \emph{DK-Torelli} (where DK stands for Dolgachev-Kapranov) if   $X$ can be reconstructed from
$T\langle X \rangle$.
 \end{definition}

 We have the following:

\begin{prop}\label{proptor1}
Let $C \subset \PP^2$ be a reduced plane curve. Then $C$ is LC-Torelli if and only if it is DK-Torelli. Therefore we just call it \emph{Torelli}. If $C$ is nonsingular then it is 
Torelli if and only if it is not of Sebastiani-Thom type.
\end{prop}

 \proof The first part is proved in \cite{Se}, Theorem 6.3. The last assertion is a special case of the main theorem of \cite{UY}. \endproof

\begin{rk}
\label{rkR(f)} If two reduced plane curves $C:f=0$ and $C':g=0$ are projectively equivalent, it is easy to see that the corresponding Jacobian rings $R(f)=S/J_f$ and $R(g)=S/J_g$ are isomorphic as graded $\C$-algebras. Under such an isomorphism,  the ideals $\wJ_f/J_f=H^0 _\mathbf{m}(R(f))$ and $\wJ_g/J_g=H^0 _\mathbf{m}(R(g))$ correspond to each other. However, these isomorphisms cannot be extended to isomorphisms of graded $S$-modules in general, as shown by the last claim in Proposition \ref{proptor1}.

\end{rk}

Let $d \ge 4$ and let $\V_{d,n,\kappa}\subset |\OO_{\PP^2}(d)|$ be the \emph{Severi variety} of plane reduced curves of degree $d$ having  $n$ nodes and $\kappa$ ordinary cusps. Then for each
$[C] \in \V_{d,n,\kappa}$ the sheaf $T\langle C \rangle$ is stable (Example \ref{exstab}(ii)) and its Chern classes $c_1,c_2$ only depend on $d,n,\kappa$ (see \ref{Chern}). Therefore we have 
a (set-theoretic, for the time being) map:
\[
\upsilon: \V_{d,n,\kappa} \longrightarrow M(2,c_1,c_2)
\]
where $M(2,c_1,c_2)$ is the moduli space of stable vector bundles of rank two and Chern classes $c_1,c_2$ on $\PP^2$.

\begin{prop}\label{severi}
The map $\upsilon$ is a morphism.
\end{prop}

\proof The proof is a straightforward consequence of the flatness of the relative first cotangent sheaf  with respect to families. For completeness we recall it. Let 
\[
\xymatrix{
\CC \ar[dr]_-\phi\ar@{^(->}[r]& \PP^2\times S \ar[d] \\
&S}
\]
be a family of curves of degree $d$ having $n$ nodes and $\kappa$ cusps, parametrized by a scheme $S$.  To this diagram one can associate the \emph{relative first cotangent sheaf}  $\mathcal{T}^1(\phi)$ which sits in an exact sequence of coherent $\OO_\CC$-modules:
\[
\xymatrix{
\OO_\CC(1)^3 \ar[r]^-\partial&\OO_\CC(d) \ar[r] & \mathcal{T}^1(\phi) \ar[r] &0}
\]
 The sheaf $\mathcal{T}^1(\phi)$ is locally generated by the partial derivatives with respect to $X_0,X_1,X_2$ of a local equation of $\CC$. It is flat over $S$ (\cite{Wa}, Lemma 3.3.8), and commutes with base change (\cite{Wa}, Lemma 3.3.6). Let $\mathcal{T}\langle \CC \rangle  = \ker(\partial)(-d)$. Twisting the above sequence by $\OO_{\PP^2\times S}(-d)$ we then obtain:
\[
\xymatrix{
0\ar[r]&\mathcal{T}\langle \CC \rangle \ar[r] & \OO_\CC(-d+1)^3 \ar[r]^-\partial&\OO_\CC \ar[r] & \mathcal{T}^1(\phi)(-d) \ar[r] &0}
\]
This exact sequence consists of coherent sheaves, flat over $S$. This last property  is a consequence of elementary properties of flatness (\cite{Har}, Prop. 9.1.A). Therefore its restrictions to the fibres of $\phi$ remain exact and therefore 
\[
\mathcal{T}\langle \CC \rangle\otimes \OO_{\CC(s)} = T\langle \CC(s) \rangle
\]
for all $s \in S$.  Therefore $\mathcal{T}\langle \CC \rangle$ defines a family of vector bundles over $S$ belonging to $M(2,c_1,c_2)$. This proves that $\upsilon$ is a morphism. \endproof

A curve $C$ belonging to $\V_{d,n,\kappa}$ is Torelli if and only if  $[C] = \upsilon^{-1}(\upsilon([C]))$. This property is clearly open in $\V_{d,n,\kappa}$.

\begin{ex}\label{extorelli1}
 From Proposition \ref{proptor1} it follows that for $n=\kappa=0$  the open set of Torelli curves coincides with the locus of nonsingular curves of degree $d$ which are not of Sebastiani-Thom type.  Note that, 
since nonsingular curves of Sebastiani-Thom type exist, the morphism $\upsilon$  has some    positive dimensional fibres  because a nonsingular curve   of Sebastiani-Thom type is linearly equivalent to infinitely many nonsingular curves having the same sheaf of logarithmic vector fields.
\end{ex}

  We have
\[
\dim(\V_{d,n,\kappa}) \ge \frac{d(d+3)}{2}-n-2\kappa
\]
and equality holds if $\kappa < 3d$ (\cite{Se1}, p. 261). 

On the other hand, if   $d=2s+1$ is odd then one computes easily that the second Chern class of 
$T\langle C\rangle\left(\frac{d-3}{2}\right)=T\langle C\rangle (s-1)$ is $c_{2,norm}=3s^2-\tau(C)$ and its first Chern class is zero. Therefore    (\cite{OSS}, p. 300)
\[\dim(M(2,0,3s^2-\tau(C)))= 4(3s^2-\tau(C))-3=12s^2-3-4n-8\kappa\]
  If $d=2s$ is even then $c_{2,norm}= c_2\left(T\langle C\rangle (s-2)\right)=3s^2-3s+1-\tau(C)$ and (\cite{OSS}, p. 317)
\[\dim(M(2,-1,3s^2-3s+1-\tau(C)))= 4(3s^2-3s+1-\tau(C))-4 = 12s^2-12s-4n-8\kappa\]

\begin{ex}\label{extorelli1}
(i) Consider the case $d=5$, i.e. $s=2$. Then the previous computations give $c_{2,norm}=12-\tau(C)$ and
\[
\dim(M(2,0,12-\tau(C)))= 45-4n-8\kappa
\]
while
\[
\dim(\V_{5,n,\kappa}) = 20-n-2\kappa
\]
and 
\[
\dim(M(2,0,12-\tau(C)))-\dim(\V_{5,n,\kappa})= 25-3n-6\kappa
\]
This implies for example that $\dim(\V_{5,10,0})  > \dim(M(2,0,2))$ and therefore that all  nodal arrangements of 5 lines are not Torelli. This is well known (see \cite{DK,jV}). It also implies that $C$ is not Torelli if $(n,\kappa)=(9,0)$ (union of an irreducible conic and three general lines).  Related examples are computed in \cite{An}, where it is shown that the nodal union of a conic with two lines is not Torelli.

(ii) Take $C$ to be the dual of a nonsingular cubic. Then $(d,n,\kappa)=(6,0,9)$ and we obtain:
$\dim(M(2,-1,1)) = 0$, while $\dim(\V_{6,0,9})=9$.  Therefore $\upsilon$ is constant and $C$ is an example of an irreducible singular curve which is not Torelli. More precisely, 
since $T\langle C\rangle(1)$ has the same Chern classes of $T_{\PP^2}(-2)$ and $M(2,-1,1)$ is irreducible \cite{Hu}, it follows that $T\langle C\rangle(1)=T_{\PP^2}(-2)$ or, equivalently, that
$T_{\PP^2}\langle C\rangle=T_{\PP^2}(-3)=\Omega^1_{\PP^2}$.
 \end{ex}

%So far these elementary computations have been valuable for discovering examples of curves that are not Torelli, but not for proving that any class of curves is Torelli. 

Using the bound on the number of nodes for irreducible curves with only nodes $(\kappa=0)$, i.e. $2n\leq (d-1)(d-2)$, we have  $\dim(M(2,c_1,c_2))\ge \dim(\V_{d,n,0})$ for any $d \geq 4$. This  induces us to conjecture that {\it the general irreducible nodal curve is Torelli}  for any $d \geq 4$ . We can prove this conjecture for irreducible curves with a small number of nodes. More precisely, we have the following result.

\begin{thm}
\label{Torthm}
Let $C:f=0$ be a degree $d \geq 4$ irreducible (resp. reducible) nodal curve in $\PP^2$, having 
$n= |\Sigma _f|>0$ nodes. If $n  < (d-1)/2$ (resp. $n  < (d-2)/2$), then $C$ is Torelli.
\end{thm}

This result will be derived from the following more technical and precise statement.

\begin{thm} \label{nodalthm} Let $C$ be a nodal curve in $\PP^2$ of degree $d\geq 4$. Let $\NN=\Sigma_f$ be the set of nodes of $C$ and consider the linear system $I_m(C)$ of curves of degree $m$ passing through the nodes in $\NN$. Assume that there is an integer $m$ such that the following holds.

\medskip

\noindent (i) $2m <d-1$ for $C$ irreducible, or $2m <d-2$ for $C$ reducible;

\medskip

\noindent (ii) the base locus of the linear system $I_m(C)$ is $0$-dimensional.

\medskip

\noindent Then the curve $C$ is Torelli.

	\end{thm}

For the proof of this result we follow essentially the same approach as in \cite{UY}, which consists of {\it  two distinct steps}: in the first step one describes in Lemma 2 the $(d-1)$st homogeneous component $J_{f,d-1}$ of the Jacobian ideal $J_f$ in terms of the sheaf
$T\langle C \rangle$, and then, in Lemma 3, one shows that the equality $J_{f,d-1}=J_{g,d-1}$ implies $f=g$, modulo a nonzero multiplicative factor, unless $f$ is of Sebastiani-Thom.

\bigskip

We consider now the first step.
Let $C:f=0$ be a reduced curve in $\PP^2$ of degree $d$, and $E:g=0$ be a (possibly nonreduced) curve in $\PP^2$ of degree $d-1$. For any $k \in \Z$, consider the exact sequence
$$ 0 \to  \OO_{\PP^2}(k-d+1) \to  \OO_{\PP^2}(k) \to  \OO_{E}(k) \to 0,$$
where the first morphism is induced by the multiplication by $g$.
Tensor this sequence by the locally free sheaf $T\langle C\rangle$ and get a new short exact sequence
$$0 \to T\langle C\rangle(k-d+1) \to T\langle C\rangle(k) \to  T\langle C\rangle (k)\otimes\OO_E\to 0.$$
The associated long exact sequence of cohomology groups looks like
$$0 \to H^0(T\langle C\rangle(k-d+1)) \to H^0(T\langle C\rangle(k)) \to 
H^0 (T\langle C\rangle (k)\otimes\OO_E) \to $$ $$ \to H^1(T\langle C\rangle(k-d+1)) \to H^1(T\langle C\rangle(k)) \to$$
Then, using the formula \eqref{logseq1}, we see that 
$$\delta_k= \dim  H^0(T\langle C\rangle(k))-\dim H^0(T\langle C\rangle(k-d+1))= 
\dim AR(f)_{k+1}-\dim AR(f)_{k-d+2}$$
 depends only on $f$ but not on $g$. Next note that the morphism
$$ H^1(T\langle C\rangle(k-d+1)) \to H^1(T\langle C\rangle(k)) $$
in the above exact sequence can be identified, using the formulas (5) and (9) in \cite{Se}  with the morphism
$$g^*_{k+1}: (\wJ_f/J_f)_{k+1} \to (\wJ_f/J_f)_{k+d}$$
induced by the multiplication by $g$. The above proves the following.
\begin{lem} \label{keylemma}
$\dim H^0 (T\langle C\rangle (k)\otimes\OO_E) = \delta_k + \dim \ker g^*_{k+1}.$
\end{lem}

Assume we are now in the situation of Theorem \ref{nodalthm}. 
Since $C$ is nodal, it follows that the saturation $\wJ _f$ of the Jacobian ideal $J_f$ coincides with the radical of $J_f$. In other words, one has $\wJ _{f,m}=I_m(C)$.
Since clearly $m<d-1$, it follows that $g^*_m$ is defined on $I_m(C)$ (considered as a vector space, not as a projective one). If $g \in J_f$, then clearly $g^*_m=0$, and hence its kernel has maximal possible dimension.

Suppose now conversely that $g^*_m=0$. By the condition (ii), it follows that there are two elements $h_1,h_2 \in I_m(C)$ having no irreducible factor in common. Since $g^*_m(h_1)=0$, it follows that $gh_1=a_1f_x+b_1f_y+c_1f_z$ for some polynomials $a_1,b_1,c_1 \in S_m$.
Similarly, we get $gh_2=a_2f_x+b_2f_y+c_2f_z$ for some polynomials $a_2,b_2,c_2 \in S_m$.
It follows that
$$(a_1h_2-a_2h_1)f_x+(b_1h_2-b_2h_1)f_x+(c_1h_2-c_2h_1)f_x=0.$$
The condition (i) implies that the only syzygy of degree $2m$ is the trivial one, i.e.
$a_1h_2=a_2h_1$, $b_1h_2=b_2h_1$ and $c_1h_2=c_2h_1$. These relations imply that $a_1,b_1,c_1$ are divisible by $h_1$, and hence $g$ is a linear combination of $f_x,f_y,f_z$.

It follows that $g \in J_{f,d-1}$ if and only if 
$$\dim H^0 (T\langle C\rangle (m-1)\otimes\OO_E) = \delta_{m-1} + \dim I_m(C),$$
i.e. the sheaf $T\langle C\rangle$ determines the homogeneous component $J_{f,d-1}$ of the Jacobian ideal $J_f$, and this completes the first step in our proof of Theorem \ref{nodalthm}. 

 Our next result, needed to complete the second step, extends Lemma 3 in \cite{UY} to a class of singular curves. Note that  in fact the hypothesis $J_f=J_g$  in  Lemma 3 in \cite{UY}  can be replaced by $J_{f,d-1}=J_{g,d-1}$, as all arguments in loc. cit. involve just linear combinations of  first order partial derivatives of some homogeneous polynomials of degree $d$.

\begin{lem}
If two irreducible distinct divisors $C$ and $D$ in $\PP^2$ of degree $d\geq 3$ and having only isolated singularities with either Milnor numbers $<(d-2)(d-1)$ or with multiplicities $<d-1$ satisfy $J_{f,d-1}=J_{g,d-1}$, then their defining equations are of Sebastiani-Thom type.
\end{lem}

\proof 

Let $C:f=0$ and $D:g=0$ be the equations of the two divisors, and assume that $f$ is irreducible and $f$ and $g$ are not proportional, as $C \ne D$. Let $\nabla h$ denote the column vector formed by the partial derivatives $h_x,h_y,h_z$ for any polynomial $h$. The equality $J_{f,d-1}=J_{g,d-1}$ implies the existence of a $3 \times 3$ constant matrix $A$ such that 
$$\nabla g =A \nabla f.$$
Let $\lambda$ be an eigenvalue of $A$ and consider the polynomial $F=g-\lambda f$. Then $\nabla F =(A-\lambda I)\nabla f$ and hence $k =\dim \langle F_x,F_y,F_z \rangle <3.$
Since $f$ and $g$ are not proportional, we have $k>0$.

\bigskip

{\bf Case 1. } $k=2$. Then by a linear coordinate change we may suppose $F_x=0$ (i.e. $F$ is a polynomial in $y,z$) and $F_y$ and $F_z$ linearly independent. The inclusion $J_F \subset J_f$ implies the existence of a $3 \times 3$ constant matrix $B$ such that 
$$\nabla F =B \nabla f.$$
Since $C$ is irreducible, it follows that $C$ is not a cone, and hence the first row in $B$ is zero.
Let $(a,b,c)$ and $(a',b',c')$ be the two other rows in $B$.

\medskip

{\bf Case 1.1 } $bc'-b'c \ne 0$. Then exactly as in the proof of Lemma 3 in \cite{UY} one shows that $f$ is of Sebastiani-Thom type.

\medskip

{\bf Case 1.2 } $bc'-b'c = 0$. The matrix $B$ has rank $k=2$, and hence we can write $f_x$ as a linear combination of $F_y$ and $F_z$, in particular $f_x$ is independent of $x$. It follows that
$f=f_0(y,z)+f_1(y,z)x$, for some homogeneous polynomials $f_0$ and $f_1$ in $y,z$.
The relation $F_y=af_x+bf_y+cf_z = af_1+b(f_{0,y}+f_{1,y}x)+c(f_{0,z}+f_{1,z}x)$ implies that
$bf_{1,y}+cf_{1,z}=0$. 

If $(b,c) \ne (0,0)$, it follows that we can make a new coordinate change involving only $y$ and $z$ such that $f_1$ becomes independent of (the new) $y$, i.e. we can take $f_1=sz^{d-1}$ for some $s \in \C$. It follows that the local equation of $C$ at the point $(1,0,0)$ is 
$sz^{d-1}+ f_0(y,z)=0$, with Milnor number at least $\mu(z^{d-1}+y^d)=(d-2)(d-1)$.

If $(b,c) = (0,0)$, then $k=2$ implies $(b',c') \ne (0,0)$ and we can repeat the same argument.

\bigskip

{\bf Case 2. } $k=1$. Then by a linear coordinate change we may suppose $F_x=0, F_y=0$ (i.e. $F$ is a polynomial in $z$) and  $F_z \ne 0$. As above we obtain a relation
$$F_z=af_x+bf_y+cf_z.$$

{\bf Case 2.1 } $c \ne 0$. Then exactly as in the proof of Lemma 3 in \cite{UY} one shows that $f$ is of Sebastiani-Thom type.

\medskip

{\bf Case 2.2 } $c = 0$. Then the relation becomes essentially $z^{d-1}=af_x+bf_y$.

Suppose first that $b=0$. Then $k=1$ implies $a \ne 0$ and hence by integration with respect to $x$ we get $af=xz^{d-1}+f_0(y,z)$. We conclude as above looking at the  local equation of $C$ at the point $(1,0,0)$. Suppose now that $b \ne 0$. If $a=0$ then we conclude as before, since the situation is now symmetric in $x,y$. 
Consider now the case when both $a$ and $b$ are nonzero. Then a linear coordinate change involving only $x,y$  brings us back to the case $a=1$ and $b=0$.

\endproof 

To complete the proof of Theorem \ref{nodalthm}, note that a Sebastiani-Thom curve in $\PP^2$ is given essentially by an equation $f_0(x,y)+z^d=0$, with $f_0$ homogeneous of degree $d$. The singular points of such a curve a given by the multiple factors of the binary form $f_0$. A factor of multiplicity $e>1$ will produce a singularity with a local equation
$u^e+v^d=0$, hence with Milnor number at least $(e-1)(d-1) \geq d-1$.
This ends the proof of Theorem \ref{nodalthm}.

To derive Theorem \ref{Torthm} from Theorem \ref{nodalthm}, consider the case when $C$ is irreducible (the other case is similar) and let $m$ be the largest integer such that $m<(d-1)/2$. If the base locus of the linear system $I_m(C)$ has positive dimension, it follows that $\dim I_m(C) \leq \dim S_{m-1}$ and hence
$$n \geq \dim S_m/I_m(C) \geq {m+2 \choose 2} - {m+1 \choose 2} =m+1.$$
Hence $m \leq n -1 < (d-1)/2-1$, a contradiction with the choice for $m$.

For curves with nodes and cusps we have the following result, by analogy to Theorem \ref{nodalthm}.

\begin{thm} \label{nodesandcuspsthm}
Let $C$ be a curve in $\PP^2$ of degree $d\geq 4$ having only nodes and cusps. Let $\NN$ be the set of nodes of $C$, $\CC$ the set of cusps and consider the linear system $I_m(C)$ of curves of degree $m$ passing through the nodes in $\NN$, the cusps in $\CC$ and having the line $T_pC$ of $C$ at $p$ as a tangent line at $p$ for any cusp $p \in \CC$. Assume that there is an integer $m$ such that
the following holds.
\medskip

\noindent (i) $2m <5d/6-2$.

\medskip

\noindent (ii) the base locus of the linear system $I_m(C)$ is $0$-dimensional.

\medskip

\noindent Then the curve $C$ is Torelli. In particular, if $C$ is a curve in $\PP^2$ of degree $d\geq 4$ having $n$ nodes and $\kappa$ cusps such that $\tau(C)=n+2\kappa \leq 5d/12-1$, then $C$ is Torelli.

	\end{thm}

\proof The proof is similar to the proof of Theorem \ref{nodalthm}, one has just to notice that the definition of the linear system $I_m(C)$ is changed in order to have again the key equality $\wJ _{f,m}=I_m(C)$. The inequality in (i) comes from the fact that in this case $\alpha_C=5/6$ as explained in Example 1.2 (ii).

\endproof

\begin{ex}  (a) Let $C$ be an irreducible curve having a unique node, say at $p=(0:0:1)$.
Then $I_1(C)=(x,y)$ satisfies the assumptions, hence $C$ is Torelli if its degree is at least $4$.
This result is sharp, since a nodal cubic is not Torelli. This follows exactly as above in Example  \ref{extorelli1}, using the inequality 
$$ \dim (M(2,0,2))=5 <\dim (V_{3,1,0})=8.$$

\medskip

(b) Let $C$ be an irreducible curve having two nodes, say at $p=(0:0:1)$ and $q=(0:1:0)$.
Then $I_2(C)=(x^2,xy,xz, yz)$ satisfies the assumptions, hence $C$ is Torelli if its degree is at least $6$.

\medskip

(c) Let $C$ be an irreducible curve having three nodes. Then there are two cases. Suppose first that the nodes are not colinear, say they are located at $p=(0:0:1)$, $q=(0:1:0)$ and $r=(1:0:0)$.
Then $I_2(C)=(xy,xz, yz)$ satisfies the assumptions, hence $C$ is Torelli if its degree is at least $6$.

When the nodes are colinear, say located at $p=(0:0:1)$, $q=(0:1:0)$ and $r=(0:1:1)$, then $I_3(C)$ contains $x^3$ and $yz(y-z)$, hence $C$ is Torelli if its degree is at least $8$.

\medskip

(d) Let $C$ be an irreducible curve having a unique cusp, say at $p=(0:0:1)$ with tangent $x=0$.
Then $I_2(C)=(x^2,xy,y^2,xz)$ satisfies the assumptions, hence $C$ is Torelli if its degree is at least $8$.

\end{ex}

%%%%%%%%%%%%%%%%%%%%%%%%%%%%%%%%%%%%%%%%%%%%%%%%%%%%%%%%%%%


\begin{thebibliography}{00}

\bibitem{An}
E. Angelini:  Logarithmic Bundles of Hypersurface Arrangements In $\PP^n$. arXiv:1304.5709.

\bibitem{ACGH}
 E. Arbarello, M. Cornalba, Ph. Griffiths, J. Harris, {\it Geometry of
Algebraic Curves I}, Springer   Grundlehren b. 267 (1984).

\bibitem{AGV}
V.I.Arnold, S.M. Gusein-Zade, A.N. Varchenko,  Singularities
  of Differentiable Maps. 2, Monographs in Math., \textbf{83},
Birkh\"auser, Basel (1988).

\bibitem{BC} R.O. Buchweitz, A. Conca: New free divisors from old. arXiv:1211.4327v1


\bibitem{D0} A. Dimca, Topics on Real and Complex Singularities, Vieweg Advanced Lecture in 
Mathematics, Friedr. Vieweg und Sohn, Braunschweig, 1987,242+xvii pp.

\bibitem{D2} A. Dimca, Syzygies of Jacobian ideals and defects of linear systems, Bull. Math. Soc. Sci. Math. Roumanie Tome 56(104) No. 2 (2013), 191- 203.

\bibitem{DS0}
A. Dimca and M. Saito,
Some remarks on limit mixed Hodge structure and spectrum, An. \c St. Univ. Ovidius Constan\c ta 22(2) (2014), 69-78.

\bibitem{DS} { A. Dimca,  M. Saito}, Graded  Koszul cohomology and spectrum of certain homogeneous polynomials, arXiv:1212.1081v3.

\bibitem{DSt0}
    A. Dimca, G. Sticlaru,  Chebyshev curves, free resolutions and rational curve arrangements,
    Math. Proc. Camb. Philos. Soc. 153, No. 3 (2012), 385-397.

\bibitem{DSt} A. Dimca, G. Sticlaru, Koszul complexes and pole order filtrations,  arXiv:1108.3976, to appear in Proc. Edinburgh Math. Soc.

%\bibitem{iD07} I. Dolgachev, Logarithmic sheaves attached to arrangements of hyperplanes, J. Math. Kyoto Univ. 47(1): 35-64 (2007).

\bibitem{DK}
 I. Dolgachev, M. Kapranov: Arrangements of hyperplanes and vector bundles on $\PP^n$. \emph{Duke
Math. J.} 71 (1993), no. 3, 633–664.

\bibitem{GMNS} M. Granger, D. Mond, A. Nieto-Reyes, and M. Schulze, Linear free divisors and the global logarithmic comparison theorem, Ann. Inst. Fourier (Grenoble) 59 (2009), no. 2, 811–850. 

%\bibitem{GLS} G.M. Greuel, C. Lossen, E. Shustin: Equisingular families of projective curves. In \emph{Global Aspects of Complex Geometry,} F. Catanese, H. Esnault, K. Hulek, A.T. Huckelberry, T. Peternell (Eds.), p. 171-209. Springer Verlag (2006).

\bibitem{Har}
R. Hartshorne: {\it Algebraic Geometry}, Springer  GTM n. 52  (1977).

\bibitem{Hu}
K. Hulek: Stable rank two vector bundles on $\PP^2$ with $c_1$ odd. \emph{Math. Annalen} 242 (1979), 241-266.

\bibitem{Ko} J. Koll\' ar: Singularities of pairs, Algebraic Geometry, Santa Cruz, 1995; Proceedings
of Symposia in Pure Math. vol. 62, AMS, 1997, pages 221-287.

\bibitem{Mac}  L. Macarro-Narv\'aez, Linearity conditions on the Jacobian ideal and logarithmic-meromorphic comparison for free divisors. Singularities I, 2008, vol. 474, p. 245-269.

\bibitem{OSS} C. Okonek, M. Schneider, H. Spindler: \emph{Vector Bundles on Complex Projective Spaces}. Progress in Math. n. 3, Birkhauser (1980).

\bibitem{OT}
P. Orlik and H. Terao,
Arrangements of Hyperplanes,
Springer-Verlag, Berlin Heidelberg New York, 1992.

\bibitem{KS} K. Saito, Einfach-elliptische Singularit\" aten, Invent. Math. 23 (1974), 289-325.

\bibitem{KS2} K. Saito: Theory of logarithmic differential forms and logarithmic vector fields, J. Fac. Sci. Univ. Tokyo Sect. IA Math. 27 (1980), no. 2, 265-291.

\bibitem{ST} A. Simis, S.O. Tohaneanu: Homology of homogeneous divisors. arXiv:1207.5862

\bibitem{Se1} E. Sernesi, {\it Deformations of Algebraic Schemes},  Springer   Grundlehren b. 334 (2006).

\bibitem{Se} E. Sernesi,  The local cohomology of the jacobian ring, Documenta Mathematica  19 (2014), 541-565. 

\bibitem{St} G. Sticlaru, Free divisors versus stability and coincidence thresholds, arXiv:1401.1843.



\bibitem{lS79}
 L. Szpiro,  \emph{Equations defining space curves,}  TATA
 Institute of Fundamental Research, Bombay (1979).

\bibitem{UY}
K. Ueda, M. Yoshinaga: Logarithmic vector fields along smooth
divisors in projective spaces, \emph{Hokkaido Math. J.}   38   (2009), 409-415.

\bibitem{jV}
J. Vall\`es, Nombre maximal d'hyperplanes instables pour un fibr\'e de Steiner, \emph{Math.
Zeit.} 233 (2000), 507-514.

\bibitem{Wa} J. Wahl: Deformations of plane curves with nodes and cusps. \emph{Amer. J. of Math.} 96 (1974), 529-577.

\bibitem{Y} M. Yoshinaga: Freeness of hyperplane arrangements and related topics,
 arXiv:1212.3523.





\end{thebibliography}
\end{document}